\documentclass{amsart}
\usepackage{amsaddr}
\usepackage{hyperref}
\usepackage{etoolbox}

\usepackage{graphicx}
\usepackage{tikz}
\usepackage{tikz-cd}
\usepackage{verbatim}
\usepackage{setspace}
\usepackage{mathrsfs}

\newtheorem{theorem}{Theorem}[section]
\newtheorem{lemma}[theorem]{Lemma}

\theoremstyle{definition}
\newtheorem{definition}[theorem]{Definition}
\newtheorem{example}[theorem]{Example}
\newtheorem{examples}[theorem]{Examples}

\theoremstyle{remark}

\theoremstyle{proposition}
\newtheorem{proposition}[theorem]{Proposition}

\theoremstyle{corollary}
\newtheorem{corollary}[theorem]{Corollary}

\numberwithin{equation}{section}

 \DeclareMathOperator{\Ext}{Ext}
\DeclareMathOperator{\Hom}{Hom} \DeclareMathOperator{\End}{End}

\DeclareMathOperator{\CE}{CE}
\DeclareMathOperator{\PE}{PE}

\DeclareMathOperator{\E}{E}

\newcommand*\xbar[1]{%
\hbox{%
\vbox{%
\hrule height 0.5pt 
\kern0.5ex
\hbox{%
\kern-0.1em
\ensuremath{#1}%
\kern-0.1em
}%
}%
}%
}
\newcommand\restr[2]{{
\left.\kern-\nulldelimiterspace 
#1 
\right|_{#2} 
}}

\begin{document}
\title{Coneat Injective Modules}

\author{Mohanad Farhan Hamid}
\address{Department of Production and Metallurgy Engineering,
\\University of Technology-Iraq
\\E-mails: {70261@uotechnology.edu.iq
\newline mohanadfhamid@yahoo.com}}

 \curraddr{}

\thanks{}





\maketitle

\section*{\bf{Abstract}}
\noindent 
A module is called \emph{coneat injective} if it is injective with respect to all coneat exact sequences. The class of such modues is enveloping and falls properly between injectives and pure injectives. Generalizations of coneat injectivity, like relative coneat injectivity and full invariance of a module in its coneat injective envelope, are studied. Using properties of such classes of modules, we characterize certain types of rings like von Neumann regular and right SF-rings. For instance, $R$ is a right SF-ring if and only if every coneat injective left $R$-module is injective.

\bigskip
\noindent \textbf{\keywordsname.} coneat submodule, coneat injective module, pure injective module, SF-ring.

\section{\bf{Introduction}}
Let $R$ be an associative ring with an identity. Unless otherwise stated, modules and $R$-modules are left unital modules over $R$. A submodule $A \subseteq B$ is called a \emph{pure} submodule (respectively, \emph{coneat} submodule) if the sequence $0 \rightarrow A \hookrightarrow B \rightarrow B/A \rightarrow 0$ is \emph{pure exact} (respectively, \emph{coneat exact}), i.e. it remains exact when tensored with any (simple) right $R$-module \cite{Cohn} (respectively, \cite{SCrivei}). In \cite{SCrivei} coneatness is studied for modules over commutative rings $R$, but here we do not assume commutativity unless stated otherwise. In \cite{ICrivei}, coneatness is called \emph{s-purity}. A class $\mathcal{F}$ of modules is called \emph{enveloping} if every module $M$ has an $\mathcal{F}$-envelope, i.e. a map $M \rightarrow E$ into a module $E$ in the class $\mathcal{F}$ such that for any other map $M \rightarrow E'$ into a module $E'$ in $\mathcal{F}$ there is a map $f: E \rightarrow E'$ that makes the diagram
$$\begin{tikzcd}
M \arrow{r}\arrow{d}{}
&E \arrow[dashed]{ld}{f}\\
E'
\end{tikzcd}$$
commutative, and such that if $E=E'$ then $f$ must be an isomorphism. Hence, envelopes are unique up to isomorphism. If the class $\mathcal{F}$ contains injectives then $\mathcal{F}$-envelopes are monomorphisms \cite[p.129]{EandJ}. In this case, we will speak of $E$, rather than the map $M \rightarrow E$, as the $\mathcal{F}$-envelope of $M$. Recall that a module is called \emph{pure injective} if it is injective with respect to all pure exact sequences. The class of pure injective modules is enveloping. For another definition of pure injective envelopes (using pure essential extensions) see \cite{FuchsSalce}. 

A module $M$ is called \emph{quasi (pure) injective} if it is injective with respect to all (pure) exact sequences having $M$ as a middle term, see \cite{JohnsonWong} and \cite{MD}. After Johnson and Wong's result \cite{JohnsonWong} that a module is quasi injective exactly when it is fully invariant in its injective envelope, a natural question arises is how much of this can be translated to quasi pure injectivity \cite[Remark 3.14]{MD}?
M. S. Abbas and the author \cite{(Flat)modules} studied modules fully invariant in their pure injective envelopes (\emph{pure quasi injective} modules) and showed that they are always quasi pure injective, but not conversely. Unlike quasi pure injectivity, the stronger concept retains some properties of quasi injectivity. For example the endomorphism ring of a pure quasi injective module is regular and self injective modulo its Jacobson radical. Another property is that the direct sum of a finite number of copies of a pure quasi injective module is always pure quasi injective \cite{(Flat)modules}.

In this paper, although we study the special case of coneat (quasi) injectivity, the method can be generalized to absorb pure (quasi) injectivity and (quasi) injectivity, as well. Here, a module is called \emph{coneat injective} if it is injective with respect to all coneat exact sequences. Every module has a coneat injective envelope containing it as a \emph{coneat essential} submodule. For the endomorphism ring and the direct sum of a finite number of copies of a coneat quasi injective module (i.e. a module that is fully invariant in its coneat injective envelope) we get analogous results to those of (pure) quasi injectives.
\noindent 
\section{\bf{Coneat Injective Modules}}
We start with the following lemma due to E. Enochs and O. Jenda \cite{EandJ}.
\begin{lemma} \label{lemmaEandJ}
\cite[p.140]{EandJ}
For every set $\mathcal{G}$ of right $R$-modules there is a unique enveloping class $\mathcal{F}$ of left $R$-modules consisting of all modules $F$ that are isomorphic to direct summands of products of copies of the left $R$-modules $G^+$ for $G \in \mathcal{G}$ and satisfying the following
\begin{enumerate}
\item $F \in \mathcal{F}$ if and only if $F$ is injective with respect to all maps $A \rightarrow B$ for which $0 \rightarrow G \otimes A \rightarrow G \otimes B$ are exact for all $G \in \mathcal{G}$.
\item A map $A \rightarrow B$ is such that $0 \rightarrow G \otimes A \rightarrow G \otimes B$ is exact for all $G \in \mathcal{G}$ if and only if all $F \in \mathcal{F}$ are injective with respect to the map $A \rightarrow B$. 
\end{enumerate}
\end{lemma}
In the above lemma, we say that the class $\mathcal{F}$ is \emph{determined} by the set $\mathcal{G}$ \cite[p.139]{EandJ}. So for example, the class of injective modules is determined by the set $\{R\}$ and the class of pure injective modules is determined by the set of modules of the form $R/I$ with $I$ a finitely generated right ideal $R$.

Now consider the set $\{R/M, M $ is a maximal right ideal of $R$ or the zero ideal$\}$.
The class $\mathcal{C}$ determined by this set consists of modules injective with respect to all coneat exact sequences and is equal to the class of left $R$-modules isomorphic to direct summands of products of copies of $G^+$ with $G$ being simple right $R$-modules.

\begin{definition} Modules in $\mathcal{C}$ above are called \emph{coneat injective} modules.
\end{definition}

From this we get the following Theorem.
\newpage
\begin{theorem} \label{first}
\
\begin{enumerate}
\item A module $M$ is coneat injective if and only if it is injective with respect to all coneat exact sequences.
\item A sequence $0 \rightarrow A \rightarrow B$ of left $R$-modules is a coneat exact sequence if and only if all coneat injective $R$-modules are injective with respect to it.
\item Every module $M$ admits a (unique up to an isomorphism) coneat injective envelope denoted $\CE(M)$. This envelope contains $M$ as a coneat submodule.
\end{enumerate}
\begin{proof}
(1) and (2) are immediate from lemma \ref{lemmaEandJ}. (3) The class of coneat injective modules is enveloping and contains the injectives as a subclass, therefore, the map $M \rightarrow \CE(M)$ is monic \cite[p.129]{EandJ}. Moreover, for any coneat injective module $C$ and any map $M \rightarrow C$ there is, by definition of envelopes, a map $\CE(M) \rightarrow C$ completing the diagram
$$\begin{tikzcd}
0 \arrow{r}& M \arrow{r}\arrow{d}{}
&\CE(M) \arrow[dashed]{ld}\\
&C
\end{tikzcd}$$
i.e. every coneat injective module is injective with respect to the sequence $0 \rightarrow M \rightarrow \CE(M)$, hence it is a coneat exact sequence by (2).
\end{proof}
\end{theorem}

\begin{corollary} \label{embed}
Every module is embedded as a coneat submodule of a coneat injective module.
\end{corollary}

The concept of coneat injectivity is weaker than injectivity and stronger than pure injectivity:
$$\text{injective} \Rightarrow \text{coneat injective} \Rightarrow \text{pure injective}$$
and it inherits some of their properties, as in the following proposition which is easy to prove.

\begin{proposition}
The class of coneat injective modules is closed for direct summands, direct products and finite direct sums.
\end{proposition}
\begin{corollary}
A module is coneat injective if and only if it is a direct summand of every module containing it as a coneat submodule.
\begin{proof} ($\Rightarrow$) Suppose $A$ is a coneat  submodule of a module $B$. If $A$ is coneat injective then the identity map of $A$ extends to a map $B \rightarrow A$ meaning that $A$ is a direct summand of $B$. ($\Leftarrow$) Conversely, embed $A$ in $\CE(A)$ as a coneat submodule so that $A$ is by assumption a direct summand of a coneat injective module and by the above Theorem it is coneat injective.
\end{proof}
\end{corollary}

Call a module $C$ \emph{coneat flat} if there are modules $A$ and $B$ with $A$ is coneat in $B$ and $C \cong B/A$. Now we have the following.
\begin{proposition} A module $M$ is coneat injective if and only if $\Ext^1 (C,M) = 0$ for all coneat flat modules $C$.
\begin{proof} Any coneat exact sequence $0 \rightarrow A \hookrightarrow B \rightarrow C \rightarrow 0$ of modules, where $C$ here must be coneat flat, gives rise to the exact sequence $\Hom(B,M) \rightarrow \Hom(A,M) \rightarrow \Ext^1 (C,M)$. Now it is clear that $\Ext^1 (C,M) = 0$ if and only if $M$ is coneat injective.
\end{proof}
\end{proposition}

\begin{corollary}
A module $C$ is coneat flat if and only if $\Ext^1 (C,M) = 0$ for every coneat injective module $M$.
\end{corollary}

\begin{proposition} The class of coneat injective modules is closed under extensions.
\begin{proof} For any short exact sequence $0 \rightarrow L \hookrightarrow M \rightarrow N \rightarrow 0$ of modules with $L$ and $N$ coneat injective and for any coneat flat module $C$, we get the exact sequence $\Ext^1(C,L) \rightarrow \Ext^1(C,M) \rightarrow \Ext^1(C,N)$. By assumption, the first and last terms are zero and therefore so is the middle one.
\end{proof}
\end{proposition}
\
\begin{examples}
\
\begin{enumerate}
\item On the ring $\mathbb{Z}$, the simple modules are $\mathbb{Z}_p$, $p$ is a prime number, and $\mathbb{Z}^+_p \cong \mathbb{Z}_p$. So all $\mathbb{Z}_p$ are coneat injective but not injective.
\item On the other hand, $\mathbb{Z}_{p^k}$, $k \geq 2$ can not be direct summands of $\prod \mathbb{Z}_p$. So they cannot be coneat injective, but they are pure injective indeed.
\item Put 
   $R=\left( {\begin{array}{cc}
   \mathbb Z_3 & 0 \\
   0 & \mathbb Z_4 \\
  \end{array} } \right)$. The only maximal ideal of $R$ is $M=\left( {\begin{array}{cc}
   \mathbb Z_3 & 0 \\
   0 & \mathbb Z_2 \\
  \end{array} } \right)$. Therefore, $R/M \cong \left( {\begin{array}{cc}
   0 & 0 \\
   0 & \mathbb Z_2 \\
  \end{array} } \right) \cong \mathbb Z_2 \cong \mathbb Z^+_2$. So, coneat injective modules are direct summands of direct products of $\mathbb Z_2$. On the other hand, $\mathbb Z^+_3 \cong \mathbb Z_3$ is pure injective but not of the above form, so it cannot be coneat injective.
\end{enumerate}
\end{examples}


Following \cite{FuchsSalce}, we call a module $B \supseteq A$ a \emph{coneat essential extension} of $A$ (and $A$ a \emph{coneat submodule of} $B$) if $A$ is coneat in $B$ and any map $\varphi : B \rightarrow X$ of $R$-modules with $\restr{\varphi}{A}$ monic and $\varphi (A)$ coneat in $\varphi (B)$ is a monomorphism. This is equivalent to saying that the only submodule $C \subseteq B$ with the property that $C \cap A = 0$ and $(A+C)/C$ is coneat in $B/C$ is the zero submodule. From this we get:

\begin{proposition} For any module $M$, the inclusion $M \hookrightarrow \CE (M)$ is coneat essential.
\begin{proof} 
Put $C = \CE(M)$ and let $\varphi$ be any map of $C$ with the property that $\varphi (M)$ is coneat in $\varphi (C)$ and $\restr {\varphi}{M}$ is monic. In the following diagram
$$\begin{tikzcd}
 M \arrow[hook]{r}\arrow{d}{\varphi}
&C \arrow{d}{\varphi}\\
\varphi (M) \arrow[hook]{r} \arrow{d}{\varphi ^{-1}}& \varphi (C) \arrow[dashed]{ldd}{\psi}\\
M \arrow[hook]{d}\\C
\end{tikzcd}$$
where existence of $\psi$ is guaranteed by coneat injectivity of $C$ and that the inclusion $\varphi (M) \hookrightarrow \varphi (C)$ is coneat, the map $\psi \varphi : C \rightarrow C$ is an automorphism since $M \rightarrow C$ is an envelope and therefore $\varphi$ must be monic.
\end{proof}
\end{proposition}

Now, following the same line of arguments in \cite{FuchsSalce} and replacing the word RD-essential by coneat essential, we get the following.

\begin{proposition} For modules $A \subseteq C$, the following are equivalent:
\begin{enumerate}
\item The module $C$ is the coneat injective envelope of $A$.
\item The module $C$ is coneat injective and is a coneat essential extension of $A$.
\item The module $C$ is the smallest coneat injective module that contains $A$ as a coneat submodule.
\item The module $C$ is a maximal coneat essential extension of $A$.
\end{enumerate}
\end{proposition}
\section{\bf{A Generalization: Relative Coneat Injectivity}}

\begin{definition}
Let $M$ and $N$ be left $R$-modules. The module $M$ is called
\begin{enumerate}
\item \emph{coneat}-$N$-\emph{injective} if for any map $f: \CE(N) \rightarrow \CE(M)$, $f(N) \subseteq M$.
\item \emph{coneat quasi injective} if $M$ is coneat-$M$-injective, i.e. if $M$ is fully invariant in its coneat injective envelope.
\item $N$-\emph{coneat injective} if for any coneat submodule $K$ of $N$, every diagram 
$$\begin{tikzcd}
0 \arrow{r}& K \arrow{r}\arrow{d}{}
&N \arrow[dashed]{ld}\\
&M
\end{tikzcd}$$
can be completed.
\item \emph{quasi coneat injective} if it is $M$-coneat injective.
\end{enumerate}
\end{definition}

The following proposition is clear.
\begin{proposition} 
\
\begin{enumerate}
\item Let $M$ and $N$ be $R$-modules. Then, $M$ is coneat-$N$-injective if and only if for any map $f: N \rightarrow \CE(M)$, $f(N) \subseteq M$.
\item Any coneat injective module is coneat-$N$-injective and any coneat-$N$-injective module is $N$-coneat injective for any module $N$.
\item A module $M$ is coneat injective if and only if it is coneat-$N$-injective if and only if it is $N$-coneat injective, for all $R$-modules $N$.
\end{enumerate}
\end{proposition}

\begin{examples}
\
\begin{enumerate}
\item In the terminology of \cite{ICrivei}, modules in the class $\mathscr{M}$ are precisely the $R$-coneat injective ones.
\item Since $\mathbb{Z}_4$ is quasi injective, it is quasi coneat injective. But as we have seen before, it is not coneat injective.
\item The module $\mathbb{Z}_3$ as a  
   $\left( {\begin{array}{cc}
   \mathbb Z_3 & 0 \\
   0 & \mathbb Z_4 \\
  \end{array} } \right)$-module is quasi injective, hence it is quasi coneat injective, but not coneat injective.
\end{enumerate}
\end{examples}

\begin{theorem} The following are equivalent for a module $M$:
\begin{enumerate}
\item The module $M$ is coneat split (i.e. every coneat submodule of $M$ is a direct summand).
\item For any module $N$, $N$ is $M$-coneat injective.
\item For any coneat submodule $N$ of $M$, $N$ is $M$-coneat-injectve.
\end{enumerate}
\begin{proof} 
(1) $\Rightarrow$ (2) $\Rightarrow$ (3) are clear. (3) $\Rightarrow$ (1) The identity map of $N$ can be extended to a map $M \rightarrow N$ and therefore $N$ is a direct summand of $M$. 
\end{proof}
\end{theorem}

\begin{theorem}  The following are equivalent for modules $M$ and $B$:
\begin{enumerate}
\item The module $M$ is $B$-coneat injective.
\item Any exact sequence $0 \rightarrow M  \overset{\alpha}{\longrightarrow} N$ splits whenever there is an $R$-map $\beta : B \rightarrow N$ with $\alpha (M) + \beta (B) = N$ and $\beta ^{-1} (\alpha (M))$ is a coneat submodule of $B$.
\end{enumerate}
\begin{proof}
(1) $\Rightarrow$ (2) Put $A = \beta ^{-1} (\alpha (M))$ and define $f: A \rightarrow M$ by $x \mapsto \alpha ^{-1} (\beta (x))$, which extends by (1) to a map $g: B \rightarrow M$. Now define $h: N \rightarrow M$ by $h(\alpha (m) + \beta (b)) = m + g(b)$.
It is easy to check that $h$ is a well-defined homomorphism with $h \circ \alpha = 1_M$.
(2) $\Rightarrow$ (1) Let $A$ be a coneat submodule of $B$ and $f: A \rightarrow M$ be a map and consider the pushout diagram
$$\begin{tikzcd}
A \arrow [hook]{r}\arrow{d}{f}
&B \arrow{d}{\beta}\\
M \arrow{r}{\alpha}& N=(M \oplus B)/W
\end{tikzcd}$$
Indeed, $\alpha (M) + \beta (B) = N$ and $\beta ^{-1} (\alpha (M)) = A$ is a coneat submodule of $B$. So, the sequence $0 \rightarrow M  \overset{\alpha}{\longrightarrow} N$ splits. Therefore, there exists an $h: N \rightarrow M$ such that $h \circ \alpha = 1_M$. So for any $a \in A$, since $\beta (a) = \alpha (f(a))$, we have $h(\beta (a)) = h(\alpha (a)) = f(a)$, i.e. $h \alpha$ extends $f$, as desired.
\end{proof}
\end{theorem}

\begin{proposition} \label{relative}
Let $M$ and $N$ be $R$-modules, then:
\begin{enumerate}
\item If $M$ is coneat-$N$-injective then $M$ is coneat-$K$-injective for all coneat submodules $K$ of $N$.
\item If $M$ is $N$-coneat injective then $M$ is $K$-coneat injective for all coneat submodules $K$ of $N$.
\item If $M$ is coneat-$N$-injective ($N$-coneat injective) then so is any direct summand $K$ of $M$.
\end{enumerate}
\begin{proof} 
\begin{enumerate} \item Any map $f: K \rightarrow \CE(M)$ extends to a $g: N \rightarrow \CE(N)$ and by assumption, $f(K) =g(K) \subseteq g(N) \subseteq M$. \item Let $A$ be a coneat injective submodule of $K$, hence it must be coneat in $N$, so that any map $A \rightarrow M$ extends to a map $N \rightarrow M$. Restrict this last map to $K$. \item is clear. \end{enumerate}

\end{proof}
\end{proposition}

\begin{corollary}
If $M_1 \oplus M_2$ is quasi coneat injective then $M_i$ is $M_j$-coneat injective for all $i, j \in \{1,2\}$.
\end{corollary}

\begin{proposition} The following are equivalent.
\begin{enumerate}
\item The direct sum of any two quasi coneat injective modules is quasi coneat injective.
\item The direct sum of any two coneat quasi injective modules is coneat quasi injective.
\item Any quasi coneat injective module is coneat injective.
\end{enumerate}
\begin{proof}
(3) $\Rightarrow$ (1) and (3) $\Rightarrow$ (2) are clear. (1) $\Rightarrow$ (3) Let $M$ be a quasi coneat injective module and let $N$ be any $R$-module, therefore $M \oplus \CE(N)$ is quasi coneat injective hence coneat injective. So $M$ is $N$-coneat injective. (2) $\Rightarrow$ (3) is similar to (1) $\Rightarrow$ (3).
\end{proof}
\end{proposition}

\begin{proposition} The module $M_1 \oplus M_2$ is coneat quasi injective if and only if $M_i$ is coneat-$M_j$-injective for all $i, j \in \{1,2\}$.
\begin{proof} ($\Rightarrow$) Clear. ($\Leftarrow$) Put $C_i = \CE(M_i)$ for $i, j \in \{1,2\}$, therefore by \cite{EandJ}, $C_1 \oplus C_2 \cong \CE(M_1 \oplus M_2)$. Let $f: M_1 \oplus M_2 \rightarrow C_1 \oplus C_2$ be any map. For all $i$ and $j$ let $f_{ij}: M_i \rightarrow C_j$ be the map $f_{ij} = \pi _j \circ \restr{f}{M_i}$. By assumption, $f_{ij} (M_i) \subseteq M_j$, so we can consider $f_{ij} : M_i \rightarrow M_j$. But $f = \left( {\begin{array}{cc}
   f_{11} & f_{21} \\
   f_{12} & f_{22} \\
  \end{array} } \right)$ applied to $\left( {\begin{array}{cc}
   m_1  \\
   m_2 \\
  \end{array} } \right) \in M_1 \oplus M_2$. Therefore, $f(M_1 \oplus M_2) \subseteq M_1 \oplus M_2$.
\end{proof}
\end{proposition}

\begin{proposition} \label{dsum}
For $R$-modules $A$ and $M = \bigoplus _{i \in I} M_i$, the module $A$ is coneat-$M$-injective if and only if $A$ is coneat-$M_i$-injective for each $i \in I$.
\begin{proof} ($\Rightarrow$) For any map $f: \bigoplus M_i \rightarrow \CE(A)$ put $f_i = \restr{f}{M_i}$. By assumption, $f(M_i) = f_i (M_i) \subseteq A$. Therefore, for every $x \in \bigoplus M_i$, $x = (x_i)_{i \in I}$ and $f(x) = \sum _{x_i \neq 0} f_i (x_i) \in A$. ($\Leftarrow$) This follows from Proposition \ref{relative}(1).
\end{proof}
\end{proposition}

Since any pure submodule is coneat, we have the following.
\begin{corollary}
If $A$ is coneat-$M_i$-injective for each $i \in I$ then $A$ is $M$-pure injective.
\end{corollary}
\begin{corollary} \label{qciqpi}
Any quasi coneat injective module is quasi pure injective.
\end{corollary}
\begin{corollary} \label{endcqi}
The endomorphism ring of a coneat quasi injective module is von Neumann regular and self injective modulo its Jacobson radical.
\begin{proof} If $A$ is a coneat quasi injective module then, by Proposition \ref{dsum}, $A$ is coneat-$A^{(I)}$-injective for every index set $I$. Hence, $A$ is $A^{(I)}$-pure injective and therefore $\End(A)$ is a cotorsion ring by \cite[Corollary 3.16]{MD}. This means, by \cite[Theorem 6]{AH}, that $\End(A)$ is regular and self injective modulo its Jacobson radical.
\end{proof}
\end{corollary}

\begin{example}
It is easy to see that the only coneat submodules of $\mathbb{Z}$ are the trivial ones. Hence, $\mathbb{Z}$ is quasi coneat injective. However, since its endomorphism ring ($\cong \mathbb{Z}$) is neither von Neumann regular nor self injective modulo its Jacobson radical, $\mathbb{Z}$ is not coneat quasi injective by corollay \ref{endcqi}.
\end{example}

The following Proposition and example show that the converse of corollary \ref{qciqpi} is not true in general.
\begin{proposition} The following are equivalent for a ring $R$:
\begin{enumerate}
\item Every coneat exact sequence of $R$-modules is pure exact.
\item All pure injective $R$-modules are coneat injective.
\item All pure injective $R$-modules are quasi coneat injective.
\item All pure injective $R$-modules are coneat quasi injective.
\item All quasi pure injective $R$-modules are quasi coneat injective.
\item All pure quasi injective $R$-modules are coneat quasi injective.
\end{enumerate}
\begin{proof} (1) $\Rightarrow$ (6) $\Rightarrow$ (5) $\Rightarrow$ (4) $\Rightarrow$ (3) are clear. (2) $\Rightarrow$ (1) Let $0 \rightarrow A \rightarrow B \rightarrow C \rightarrow 0$ be a coneat exact sequence of $R$-modules and therefore, by (2), every pure injective $R$-module is injective with respect to it. So, by \cite[p.290]{Wis}, the sequence is pure exact. (3) $\Rightarrow$ (2) Let $M$ be a pure injective $R$-module and therefore, so is $M \oplus \PE(N)$ for any $R$-module $N$. Hence $M \oplus \PE(N)$ is quasi coneat injective, by (3). So $M$ is $\PE(N)$-coneat injective. But $N$ is a pure, hence coneat, submodule of $\PE(N)$. So $M$ is $N$-coneat injective for all $R$-modules $N$.
\end{proof}
\end{proposition}

Applying the above propsition, since $\mathbb{Z}_3$ as a $\left( {\begin{array}{cc}
   \mathbb Z_3 & 0 \\
   0 & \mathbb Z_4 \\
  \end{array} } \right)$-module is pure injective but not coneat injective, we see that the ring $\left( {\begin{array}{cc}
   \mathbb Z_3 & 0 \\
   0 & \mathbb Z_4 \\
  \end{array} } \right)$ has a coneat exact sequence which is not pure exact.

Recall that a ring $R$ is called a \emph{right SF-ring} \cite{Ram} if each simple right $R$-module is flat. This is equivalent to saying that the character modules of simple right $R$-modules are injective left $R$-modules which is equivalent to saying that every coneat injective left $R$-module is injective. Therefore, every von Neumann regular ring is an SF-ring. Conversely, commutative SF-rings are von Neumann regular since so is the center of a right SF-ring \cite{Ram}.

\begin{theorem} \label{SF}
Consider the following statements for a ring $R$:
\begin{enumerate}
\item The ring $R$ is a right SF-ring.
\item Every coneat injective left $R$-module is injective.
\item Every coneat injective left $R$-module is quasi injective.
\item Every coneat quasi injective left $R$-module is quasi injective.
\item Every quasi coneat injective left $R$-module is quasi injective.
\item Every coneat injective left $R$-module is absolutely coneat (= coneat in every module contatining it as a submodule \cite{SCrivei}).
\item All left $R$-modules are absolutely coneat.
\item Every left $R$-module $M$ is coneat regular (i.e. every submodule of $M$ is coneat).
\item All exact sequences of left $R$-modules are coneat exact.
\item The ring $R$ is coneat regular.
\item All left $R$-modules in the class $\mathscr{M}$ are injective.
\item The ring $R$ is von Neumann regular.
\end{enumerate}

Then \emph{(1) -- (11)} are equivalent and \emph{(12)} $\Rightarrow$ \emph{(1)}. If $R$ is commutative then all the statements are equivalent.
\begin{proof}
(12) $\Rightarrow$ (1) $\Leftrightarrow$ (2) and, in the commutative case, (1) $\Rightarrow$ (12) are given in the above paragraph. (2) $\Rightarrow$ (6), (5) $\Rightarrow$ (4) $\Rightarrow$ (3) and (7) $\Leftrightarrow$ (8) $\Leftrightarrow$ (9) $\Rightarrow$ (10) are clear. (3) $\Rightarrow$ (2) Let $M$ be a coneat injective $R$-module. The coneat injective module $M  \oplus \CE (R)$ is quasi injective, by assumption. It is easy to see that $M$ is injective with respect to $R$ (i.e. injective). (6) $\Rightarrow$ (2) Any coneat injective module is absolutely coneat and hence a direct summand of every module containing it. (2) $\Rightarrow$ (7) Any $R$-module $M$ is coneat in the injective $R$-module $\CE(M)$ and therefore, $M$ must be absolutely coneat. (8) $\Rightarrow$ (5) Any module that is both coneat regular and coneat quasi injective must be quasi injective. (10) $\Rightarrow$ (11) The $R$-modules in the class $\mathscr{M}$ are precisely the ones injective with respect to all coneat exact sequences having $R$ as a middle term. If $R$ is coneat regular then such modules must be injective by Baer condition. (11) $\Rightarrow$ (2) Coneat injective $R$-modules are in the class $\mathscr{M}$. 
\end{proof}
\end{theorem}
\begin{theorem} \label{cssimple}
The following are equivalent for a ring $R$:
\begin{enumerate}
\item Every coneat exact sequence of $R$-modules splits.
\item All $R$-modules are coneat injective.
\item All $R$-modules are coneat quasi injective.
\item All $R$-modules are quasi coneat injective.
\end{enumerate}
\begin{proof}
(1) $\Rightarrow$ (2) $\Rightarrow$ (3) $\Rightarrow$ (4) are immediate. (4) $\Rightarrow$ (2) For any $R$-modules $M$ and $N$, the module $M  \oplus \CE (N)$ is quasi coneat injective, by assumption. Therefore, $M$ is $\CE (N)$-coneat injective and hence it is $N$-coneat injective for all $N$. So $M$ must be coneat injective, as desired. (2) $\Rightarrow$ (1) Let $0 \rightarrow A \rightarrow B \rightarrow C \rightarrow 0$ be a coneat exact sequence. Since $A$ is coneat injective, the identity map of $A$ extends to a map $B \rightarrow A$ so that the sequence splits.
\end{proof}
\end{theorem}

Combining Theorems \ref{SF} and \ref{cssimple}, we get the following characterization of semisimple artinian rings.
\begin{theorem} The necessary and sufficient condition for a ring $R$ to be semisimple artinian is that $R$ is a right SF-ring and every coneat exact sequence of left $R$-modules splits.
\end{theorem}
\bibliographystyle{amsplain}

\end{document}